\documentclass[11pt]{amsproc}
\usepackage[dvips]{graphicx,color}
\usepackage{amssymb,amscd,latexsym,epsfig}
\usepackage[mathscr]{euscript}

\DeclareFontFamily{OT1}{rsfs}{}
\DeclareFontShape{OT1}{rsfs}{n}{it}{<-> rsfs10}{}
\DeclareMathAlphabet{\curly}{OT1}{rsfs}{n}{it}

\newcommand\C{\mathbb C}
\newcommand\R{\mathbb R}
\newcommand\Z{\mathbb Z}
\newcommand\Pee{\mathbb P^1}
\newcommand\Pt{\mathbb P^{\,3}}
\newcommand\Pf{\mathbb P^{\,4}}
\newcommand\OO{\mathscr O}
\newcommand\into{\hookrightarrow}

\newcommand\res{\arrowvert_}
\newcommand\To{\longrightarrow}
\newcommand\Ree{\mathrm{Re\,}}

\makeatletter \@addtoreset{equation}{section} \makeatother

\newtheorem{Defn}[equation]{Definition}

\newtheorem{Theorem}[equation]{Theorem}
\newtheorem{Lemma}[equation]{Lemma}

\newtheorem{Prop}[equation]{Proposition}
\newenvironment{Proof}{\noindent\emph{Proof.}}{\hfill$\square$\\}
\newenvironment{Remark}{\noindent\textbf{Remark}.}{\\}
\newenvironment{Remarks}{\noindent\textbf{Remarks}.}{\\}

\title{\textbf{Symplectic conifold transitions}}
\author{I. Smith, R. P. Thomas and S.-T. Yau}
\date{}

\begin{document}

\begin{abstract} \noindent
We introduce a symplectic surgery in six dimensions
which collapses Lagrangian three-spheres and replaces
them by symplectic two-spheres. Under mirror symmetry it
corresponds to an operation on complex 3-folds studied by Clemens,
Friedman and Tian.  We describe several examples which show that there are
either many more Calabi-Yau manifolds (e.g. rigid ones) than previously
thought or there exist ``symplectic
Calabi-Yaus'' -- non-K\"ahler symplectic 6-folds with $c_1=0$. The analogous
surgery in four dimensions, with a generalisation to ADE-trees of
Lagrangians, implies that the canonical class of a
minimal complex surface contains symplectic forms if and only if it has
positive square.
\end{abstract}

\maketitle

%%%%%%%%%%%%%%%%%%%%%%%%%%%%%%%%%%%%%%%%%%%%%%%%%%%%%%%%%%%%%%%%%%%%%%%%%%%

\section{Introduction}

A 3-fold ordinary double point, or \emph{node} as we shall call it,
is a complex singularity analytically equivalent to
$$
\{xy=zw\}\subset\C^4.
$$
By taking the graph of the rational function $x/z=w/y$ from a
neighbourhood of the singularity to $\Pee$ we get a \emph{small
resolution} of the node; a smooth resolution with exceptional
set $\Pee$. This is because away from the origin at least one
of $x/z$ or $w/y$ is uniquely
defined on $\{xy=zw\}$, so the graph is isomorphic to the domain away
from the origin, and replaces the origin with the whole $\Pee$.
Similarly using the function $x/w=z/y$ gives
another small resolution, the \emph{flop} of the first.
Alternatively, smoothing the node, $\{ xy-zw=\epsilon \}$,
yields a 3-sphere vanishing cycle (described below).

So given a node on a K{\"a}hler 3-fold
one may either try to smooth it (producing a Lagrangian $S^3$
vanishing cycle) or resolve it (producing
a holomorphic $\Pee$).  Passing from one desingularisation to the
other is called
a \emph{conifold transition} in the physics literature.  There is a
\emph{natural complex structure on the
resolution}, but \emph{not} a natural K{\"a}hler structure.  Indeed, locally
there is an obvious parameter for any K{\"a}hler form, given by the
symplectic area of the resolving sphere, and the existence of such
choices means that there are also obstructions to patching together the
local choice of symplectic form to a global form.  Hence the resolution
may not admit a symplectic structure compatible with the complex structure.

Conversely the smoothing $\{xy-zw=\epsilon\}$ \emph{is
in a natural way symplectic}, but \emph{not} naturally complex.
At least in the Calabi-Yau setting, it is the volume of the vanishing
cycle, computed via integrating 
the holomorphic 3-form, or equivalently $\epsilon$, which defines a
local parameter for the choice of complex structure, and there is an
obstruction to patching these choices to give a global complex structure.
However, on any complex smoothing there is a natural compatible
symplectic structure, and even if complex smoothings do not exist
\cite{Fr}, this ``symplectic smoothing'' does (Theorem \ref{smooth} below).
As explained in Section 2, the elementary but
fundamental fact is that there is a symplectomorphism between 
$\{ xy=zw \} \backslash \{ 0 \}$ and $T^* S^3 \backslash \{
\mathrm{Zero \ section} \}$,
equipped with their standard symplectic structures.

In the Calabi-Yau case, there are necessary and sufficient conditions
\cite{Fr}, \cite{Ti} for the existence of a complex smoothing; in the
K{\"a}hler case these can be interpreted
as saying that the ``symplectic smoothing''
admits a compatible complex structure if and only if the
conditions of Friedman-Tian are satisfied.

Under mirror symmetry, the mirrors of Calabi-Yau manifolds with nodes
are usually also Calabi-Yaus with nodes, and the smoothing and resolution
processes get swapped \cite{Mo}. So there should be a
criterion mirror to that of Friedman-Tian giving necessary and
sufficient conditions for the resolution of a symplectic manifold with
nodes to admit a symplectic structure.
We give such a result here (Theorem \ref{resol}), giving a new way
to produce symplectic manifolds via the surgery of replacing
Lagrangian $S^3$s with symplectic $\Pee$s (and preserving $c_1$
in the process). 

We give a number of examples of conifold transitions that preserve
a symplectic structure. In the
complex setting, Lu and Tian
\cite{LuT} produce a complex structure
(non-K{\"a}hler, with trivial canonical bundle) on
$(S^3\times S^3)^{\#(n\geq 2)}$, and mirror to this we produce
``symplectic Calabi-Yaus" (symplectic manifolds with $c_1=0$) with Betti
numbers $b_3 =2, 2\leq b_2 \leq
25$.  One can think of this as mirroring Reid's fantasy \cite{Re}.
However, it is not clear if one can go further than this, i.e. if there can be
symplectic structures with vanishing first Chern class on manifolds with
$b_3=0$.

To set this in context, note that for simply-connected 4-manifolds $X$,
symplectic geometry is very similar to K{\"a}hler geometry when either
$c_1(X)>0$ (i.e. when $c_1(X)$ can be represented by a symplectic
form) or $c_1(X)=0$. In the first case $X$ can
be shown to be Fano, and in the second it follows from results 
of Morgan and Szabo \cite{MoSz} that $X$ is homeomorphic to the
K3 surface. Beyond this, symplectic and K{\"a}hler geometry diverge,
one reason being the existence of symplectic
surgeries -- fibre connect sums -- which are non-K{\"a}hler.  (The
analogue of our surgery in four dimensions, and its cousins, are
themselves interesting, as we point out in \ref{symplectic}.) We
expect such a divergence for $c_1=0$ in 6-dimensions, where conifold
transitions provide a symplectic
surgery preserving $c_1$.  However, finding and studying 
Lagrangian $S^3$s and their configurations -- the geometric input for
such a surgery --  is much harder than finding holomorphic $\Pee$s, and
so it can be hard to find examples for which the surgery gives a non-K{\"a}hler
result. In particular, controlling the
intersections of the Lagrangians can be highly non-trivial -- for instance,
it may be that there are 
subtle, Floer-theoretic obstructions to obtaining disjoint families of
Lagrangian spheres spanning $b_3/2$-dimensional subspaces of $H_3$, preventing
us from using our surgery to produce symplectic 6-folds with vanishing
$c_1,\,b_3$ and $\pi_1$ (which would be necessarily non-K{\"a}hler by
(\ref{notkahler})). For instance, Donaldson has asked \cite{Do} if all
Lagrangian spheres in algebraic varieties arise as vanishing cycles for
complex degenerations, and one can check that a positive answer to this
question would exclude the existence of essential Lagrangian spheres
in Calabi-Yau 3-folds with $b_3=2$.
 
\vspace{0.2cm}

\noindent \textbf{Acknowledgements.} We would like to thank 
Denis Auroux, Philip Candelas, Xenia de la Ossa and Gang Tian for useful
conversations.
The authors are supported by an EC Marie-Curie fellowship
HPMF-CT--2000-01013, a Royal Society university research
fellowship, and DOE and NSF grants DE-FG02-88ER35065 and DMS-9803347NSF.

%%%%%%%%%%%%%%%%%%%%%%%%%%%%%%%%%%%%%%%%%%%%%%%%%%%%%%%%%%%%%%%%%%%%%%%%%%

\section{Smoothings and resolutions}

To describe the symplectic
versions of degeneration and resolution, we will begin with some local
facts about nodes.  Fix once and for all a complex parametrisation $W
= \{\sum z_i^2 = 0\}$ of a node.  This complex variety has three
resolutions of singularities relevant to the discussion.  First, one
can blow up the singular point to obtain a variety $W_b$, which can be
described as follows. Blowing up the origin of $\C^4$
introduces an exceptional divisor $\Pt$; the closure of 
$W \backslash \{0\}$ inside the blow-up of $\C^4$ at $0$
meets this $\Pt$ in the quadric surface $E = \Pee \times \Pee
\subset \Pt$ given by the defining equation $\{\sum z_i^2 = 0\}$. This
closure $W_b$ is smooth and the normal bundle to $E\subset W_b$ is 
$\OO(-1, -1)$.  There are also the two \emph{small} resolutions $W^\pm$
mentioned in the Introduction given by blowing down either of the two
rulings $E\to\Pee$ of $E$. Thus each of $W^\pm$ has exceptional locus a
$\Pee$ over $0 \in W$ (with normal bundle $\OO(-1) \oplus \OO(-1)$), and
blowing up the $\Pee\subset W^\pm$ gives back $W_b$.
The projection maps of the resolutions define canonical isomorphisms
$(W_b\backslash E) \cong (W\backslash \{0\}) \cong (W^{\pm} \backslash \Pee)$.

It will be important for us to fix models of $W^{\pm}$ so that the two are
distinct and not interchangeable. Changing the coordinates in the
introduction by
$$x \mapsto z_1+iz_2, \, y \mapsto z_1 - iz_2, \, z\mapsto
-z_3-iz_4, \, w\mapsto z_3-iz_4$$
takes $\{xy=zw\}$ to $\{\sum z_i^2=0\}$, so we fix $W^\pm$
to be defined via the graphs of the following rational maps
$$W^+: \ \frac{z_1+iz_2}{z_3+iz_4} \ = \ -\frac{z_3-iz_4}{z_1-iz_2}~;
\qquad W^-: \ \frac{z_1+iz_2}{z_3-iz_4} \ = \
-\frac{z_3+iz_4}{z_1-iz_2}.$$
In particular changing the choice of coordinates on $W$ by $z_4\mapsto
-z_4$ (which preserves $\sum z_i^2$) swaps $W^\pm$.

Now we relate these resolutions to the cotangent bundle of the
three-sphere.  For the standard oriented $S^3 \subset \R^4$ we
can fix coordinates on $T^*S^3$ as follows:
$$
T^* S^3 \ = \ \{ (u,v) \in \R^4 \times \R^4 \ | \ |u| = 1, \langle
u,v \rangle = 0 \}.
$$
The key local fact is that there is a symplectomorphism
\begin{equation} \label{localmodel}
(W_b\backslash E\ \cong\ W^{\pm} \backslash \Pee\ \cong\ ) \quad\quad
(W \backslash \{0\}, \omega_{\C^4}) \ \stackrel{\phi\,}{\longrightarrow}
\ (T^*S^3\backslash \{v=0\}, d(vdu))
\end{equation}
which can be given explicitly in coordinates via the map
$$
(z_j = x_j + iy_j)_{1 \leq j \leq 4} \ \mapsto \ (x_j / |x|, -|x|
y_j)_{1 \leq j \leq 4}.
$$ 
(Here $T^*S^3\backslash \{v=0\}$ is the cotangent bundle of $S^3$ minus
its zero-section, and $|x| = (\sum x_i ^2)^{1/2}$ is the norm of the real
vector which is the real part of $z$.  One computes directly that
$\phi ^* (\sum_j dv_j \wedge du_j) = (i/2) \sum_j dz_j \wedge
d\overline{z}_j$.  For more discussion of this, see
\cite{Se}.)

From this, we have the following important observation.  The
holomorphic, hence orientation-preserving automorphism of $\C^4$ given by
$z_4 \mapsto -z_4$ (that preserves $W$ and interchanges the two
small resolutions $W^{\pm}$) acts on the real slice $\R^4 \subset \C^4$
by reflection in a
hyperplane, and in particular induces an orientation-reversing
diffeomorphism of $\{|u|=1\} \cong S^3 \subset \R^4$.  In other words,
\emph{flopping the} $\Pee \Leftrightarrow$ \emph{changing orientation
on the} $S^3$.  A more thorough
description of this, involving the relevant homogeneous spaces, moment
maps, and the topology of the surgery we are performing, is in
the Appendix.

\vspace{0.2cm}

We now globalise the transition from a smoothing to a resolution
and back again.  To obtain smooth surgeries well-defined up to
diffeomorphism it will be
important to control the choices involved; this same information will
later give us control on the symplectic structures via Moser's theorem.

Thus, we can make the above discussion relevant to more general
symplectic manifolds by recalling \cite{We} that if $L \subset X$ is
a Lagrangian three-sphere in a symplectic six-manifold, then a
neighbourhood of $L$ in $X$ is symplectomorphic to a neighbourhood of
the zero-section in $T^*L$, equipped with its canonical symplectic
structure.  (An explicit symplectomorphism can be defined by a choice
of compatible almost complex structure on $X$, and these form a 
contractible space \cite{McD-S}.)  Such a Lagrangian is not
canonically oriented. However, if we 
pick an orientation, then there is a unique orientation-preserving
diffeomorphism to $S^3 \subset \R^4$ up to homotopy, since
$\mathrm{Diff}_+ (S^3) \simeq SO(4)$ by a famous result of Hatcher \cite{H}.
This induces a symplectomorphism between a neighbourhood of $L$ and
$T^*S^3$.

Similarly, given a symplectic two-sphere $C \subset X$ with normal
bundle having first Chern class $-2$, Weinstein's neighbourhood theorem
\cite{We} implies that a neighbourhood of $C$ in $X$ is
symplectomorphic to a neighbourhood of the zero-section in the bundle
$\OO(-1) \oplus \OO(-1) \rightarrow \Pee$ equipped with a K{\"a}hler
form giving the zero-section the same area as $C$.  To define
such a symplectomorphism, following (\cite{McD-S}, p.94-5), choose
an almost complex structure taming $\omega$ and making $C$
$J$-holomorphic, together with a sufficiently small positive number
$\varepsilon$ so that the exponential map of the metric defined by
$\omega$ and $J$ is injective on the $\varepsilon$-disc bundle inside
the normal bundle of $C$. Both $C$ and
$\Pee$ are canonically oriented; also fix an orientation-preserving
diffeomorphism $C \rightarrow \Pee$ and a lift of this to a linear
isomorphism of complex normal bundles $\nu_{C/X} \rightarrow \OO(-1) \oplus
\OO(-1)$.  Given all this data, \cite{McD-S} gives an explicit
symplectomorphism from a small neighbourhood of $C$ to a small
neighbourhood of $\Pee$ inside the $\OO(-1)\oplus \OO(-1)$ bundle.
Each of the sets of required choices is connected: the space of tamed almost
complex structures is contractible, the choice of diffeomorphism belongs to 
$\mathrm{Diff}_+(S^2) \simeq SO(3)$, and any two lifts to bundle maps
differ by a change of framing of $\nu_{C/X}$; such framings are
parametrised by $\pi_2(SO(4)) = \{1\}$.  Hence, up to homotopy, there
is a unique symplectomorphism from a neighbourhood of $C$ in $X$ to a
neighbourhood of the exceptional curve $\Pee$ in either small
resolution $W^{\pm}$ of a node.  Here we have fixed complex
co-ordinates  with respect to which our surgery is canonically defined.
It follows that the operations
defined below yield manifolds well-defined up to isotopy (and in
particular diffeomorphism):

\begin{Defn}
Let $X$ be a symplectic six-manifold and $L \subset X$ a Lagrangian
three-sphere.  By a \emph{conifold transition} of $X$ in $L$, we mean
the \emph{smooth} manifold $Y=(X \backslash L) \cup_{\phi^{-1}} W^{\pm}$,
where we fix a  diffeomorphism $L \rightarrow S^3$ and hence
parametrise a neighbourhood of $L$ by a neighbourhood of the
zero-section in $T^*S^3$; then $\phi^{-1}$ is the (restriction to this
neighbourhood of the) diffeomorphism (\ref{localmodel}) $T^*S^3
\backslash \{v=0\} \rightarrow W^{\pm} \backslash \Pee$.

Let $Y$ be a symplectic manifold and $C \subset X$  a symplectic
two-sphere whose normal bundle has Chern class $-2$.  The
\emph{reverse conifold transition} of $Y$ in $C$ is the smooth
manifold $X=(Y\backslash C) \cup_{\phi} (T^*S^3)$, where $\phi$ is the
(restriction to suitable neighbourhoods of the)
diffeomorphism of (\ref{localmodel}) composed with a diffeomorphism of
$\nu_{C/Y}$ and $W\backslash \{0\}$ as above. 
\end{Defn}

By the preceeding discussion, then, up to homotopy there
are two $\Z/2\Z$ 
choices in the conifold transition -- we orient $L$, and we choose a
small resolution.  Swapping both choices gives back the same smooth
manifold since, as we have seen, changing orientation on $S^3$
interchanges the factors of $\Pee \times \Pee$ and swaps the small
resolutions, and vice-versa.  Hence, up to diffeomorphism, there are 
exactly two distinct conifold transitions. The
reverse conifold transition is uniquely defined as a smooth manifold,
but the obvious embedded three-sphere is \emph{not}
canonically oriented, since changing the gluing map via the 
automorphism $z_4 \mapsto -z_4$ on $W\backslash \{0\}$ changes the
orientation on the $S^3$.

To understand how symplectic structures interact with conifold transitions
it will be convenient to use the intermediate space with a node, and a
model symplectic structure at the node. So we make
the following definition, in $n$ complex dimensions for now; later $n$
will be 3. We say a continuous map $\phi: \C^{n+1} \rightarrow
\C^{n+1}$ is \emph{admissible} if it is smooth away from the origin,
$C^1$ at the origin with $d\phi_0\in\,$Sp\,$(2n+2,\R)$, and setwise fixes
the quadric $W= \{ \sum z_i ^2 = 0 \}$. 

\begin{Defn} \label{admiss}
A \emph{conifold} is a topological space $X$ covered by an atlas of
charts $\{(U_i, \phi_i)\}_{i\in I}$ of the following two types: either
$\phi_i~: U_i \rightarrow D^{2n}$ is a homeomorphism onto an open
disc in $\R^{2n}$ or $\phi_j~: U_j \rightarrow W \cap D^{2n+2}$ is
a homeomorphism onto the intersection of an open disc in $\C^{n+1}$
with the quadric $W = \{ \sum_{i=1}^{n+1} z_i ^2 = 0 \} \subset \C^{n+1}$.
In the latter case, the point $P = \phi_j^{-1}(0)$ is
called a \emph{node} of $X$.

Moreover, the transition maps $\phi_{ij} = \phi_i \circ
\phi_j^{-1}$ must be $C^\infty$ away from nodes, and if $P
\in U_i \cap U_j$ is a node then there must be an open subset $V
\subset \C^{n+1}$ containing $0$ such that
$\phi_{ij}\res{V\cap W}$ coincides with the restriction of an admissible
homeomorphism.
\end{Defn}

We will call such charts \emph{smooth admissible coordinates}.

\begin{Defn} \label{symponnode}
A \emph{symplectic structure} on a conifold $X$ is a smooth closed
non-degenerate two-form $\omega_X$ on $X \backslash
\{\mathrm{Nodes}\}$ which, in any set of admissible coordinates around
each node, 
co-incides with the restriction of a closed two-form on $\C^{n+1}$
which is smooth away from $0 \in \C^{n+1}$, and continuous and 
equal to the standard K{\"a}hler 2-form at the origin. 

Two such closed forms $\omega_i$ define \emph{equivalent} symplectic structures
on $X$ if there exists an admissible homeomorphism $\phi$ of $X$ such that
$\phi^*\omega_1\equiv\omega_2$ on $X \backslash \{\mathrm{Nodes}\}$.
\end{Defn}

We will call such an $(X,\omega_X)$ a symplectic conifold.  (Observe
that the class of two-forms we consider on $\C^{n+1}$ is preserved by
admissible homeomorphisms, and using such homeomorphisms necessitates
allowing forms which are only continuous at zero.)
We have a ``Darboux theorem'', asserting that locally 
the symplectic structure is unique near any node.  This would not be
possible without some pointwise information at the node, for instance,
consider the 2-dimensional case: the 1-fold node is
two symplectic two-planes meeting transversely at a point; writing one as a graph $f: \R^2 \rightarrow \R^2$ over the
symplectic orthogonal complement of the other at the node, the trace and
determinant of $df$ are local symplectic invariants.

\begin{Prop} \label{nodedata}
Let $X$ be a symplectic conifold and let $P \in X$ be a node. There is
some neighbourhood $U$ of $P$ with admissible coordinates (\ref{admiss})
in which $\omega_X$ is equivalent to the restriction of 
$\omega_{\C^{n+1}}$ to a neighbourhood of $0 \in W$.
\end{Prop}

\begin{Proof}
Fixing an admissible chart (\ref{admiss}), we may assume we are
working on a neighbourhood of the origin in the standard node $W$,
with a non-standard symplectic structure $\omega$ defined on
a ball near the origin in $\C^{n+1}$ and co-inciding with the standard
structure $\omega_{\C^{n+1}}$ at the origin.  Recall the usual proof
of the Darboux 
theorem.  On the ball around $0 \in \C^{n+1}$ we can choose a one-form
$\sigma$ such that $d\sigma = \omega - \omega_{\C^{n+1}}$, and without loss
of generality we can suppose $\sigma$ vanishes (to order two) at the
origin.  We then define a
family of vector fields $X_t$, all vanishing continuously at the origin, via
\begin{equation} \label{uselater}
\sigma + \iota_{X_t} \omega_t = 0
\end{equation}
where $\omega_t = \omega_{\C^{n+1}} + td\sigma$.  The flow of this
family of vector fields yields a family of diffeomorphisms $\{f_t\}$
of $\C^{n+1}\backslash\{0\}$, extending as $C^1$-maps over $0$ which
fix the origin, such that $f_1$ pulls back $\omega_{\C^{n+1}}$ to $\omega$.  

Note that, at least in a small enough ball around the origin, the
linear family of symplectic forms $\omega_t$ will all have
non-degenerate restriction to $W$ (e.g. all the forms tame the
integrable complex structure $J$ in some small ball and $W$ is
$J$-holomorphic). For the same reason, in some ball the forms will all
be symplectic on the 
fibres $W_t$ (where $W_0=W$) of the map $\pi: z \mapsto \sum z_i^2$.
We may
therefore define the ``horizontal projections'' $\{H_t\}$ of the
vector fields $\{X_t\}$ as follows.  For every $0\leq t\leq 1$ and $0\neq z \in
\C^{n+1}$ there is a real rank two subbundle of $T_z\C^{n+1}$ which is
the $\omega_t$-symplectic orthogonal complement to the tangent bundle
of the fibre of $\pi$ through $z$.  Let $H_t$ denote the projection of
$X_t$ to this real rank two subbundle.  This is certainly smooth as a
vector field on $\C^{n+1}\backslash\{0\}$, and we claim that it has a
continuous extension over the origin.  

For $t=0$, this follows by a direct computation.  In this case,
$\omega_0 = \omega_{\C^{n+1}}$ is the standard symplectic structure,
and the symplectic orthogonal complement to $T_z(\pi^{-1}(\pi(z)))$ is
generated by the complex conjugate vector, i.e. is $\C \langle
\overline{z} \rangle$.  It follows that the vector field $H_t$ is
given by $H_t(z) = \alpha(z)\overline{z}$ where the function $\alpha$
is defined by the identity
$$\omega_0(X_0(z)-\alpha(z)\overline{z},\overline{z}) \ = \ \langle
X_0(z)-\alpha(z)\overline{z}, i\overline{z}  \rangle  \ = \ 0 $$
with $\langle \cdot, \cdot \rangle$ the usual metric on
$\C^{n+1}$. Although the function 
$\alpha$ may not be continuous, as not every function
vanishing to order two can be written as $|z|^2$ times a continuous
function, the vector field 
$H_0(z) = \alpha(z)\overline{z}$ does vanish continuously at the origin
since $\langle X_0(z), \overline{z}\rangle$ vanishes to order two.
The vector fields $H_t$ will also vanish continuously at the origin,
since these can be determined explicitly from $H_0$ by functions of
the global changes of co-ordinates given by the family of maps $\{f_t\}$.  

Integrating up, let $\{F_t\}$ denote the flow of the vector fields
$V_t = X_t-H_t$, well-defined in a small ball around the origin.
Since the vector fields are tangent to the fibres of
$\pi$ and $C^0$ at the origin, the maps $F_t$ will all be admissible.
Given the definition of equivalence of symplectic forms on conifolds,
it is enough to prove that $F_1^* \omega$ and $\omega_{\C^{n+1}}$ have
the same restriction to the open quadric $W\backslash\{0\}$.  Using
the defining equation $\frac{d}{dt} F_t = V_t \circ F_t$ we find that 
$$\frac{d}{dt} F_t^* \omega_t \ = \ F_t^* \left( \frac{d\omega_t}{dt}
  + d\iota_{V_t}\omega_t \right);$$
combining this with $d\omega_t/dt = d\sigma$ and \ref{uselater} we have 
$$\frac{d}{dt} \left( (F_t^* \omega_t)|_{W\backslash\{0\}} \right) \ =
\ F_t^* \left( \iota_{H_t}\omega_t)_{W\backslash\{0\}} \right) \ = \
0$$
since $\omega_t(H_t, dF_t(u)) = 0$ for any $u\in \ker(d\pi)$.  We
therefore have that $F_t^* \omega_t$ is constant on restriction to the
quadric, and in particular $(F_1^* \omega)|_{W\backslash\{0\}} =
\omega_{\C^{n+1}}|_{W\backslash\{0\}}$ with $F_1$ admissible, as
required.
\end{Proof}

It is very likely that, as for smooth subvarieties $Z$, one can in
fact find a smooth change of co-ordinates taking $\omega$ to
$\omega_{\C^{n+1}}$ even in the case where there is an isolated singular
point.  This would require a more substantial analysis.
Despite the loss of regularity, the above implies that \emph{any
symplectic form near the node which is standard on the node,
is symplectomorphic, in a punctured
neighbourhood, to the restriction of the form
$\omega_{\C^{n+1}}$ to the punctured quadric $W \backslash \{0\}$}.
This is all we shall require in the sequel.  For later, note also that
the proof 
has an obvious extension to symplectic structures on manifolds with
isolated singular points modelled on other singularities, for instance
on ADE singularities.

Now, mirror to the fact that a small resolution of a node on a complex
variety is again naturally complex, we can show that the smoothing is
naturally a symplectic operation.  The proof shows that, although we
refer to a ``symplectic smoothing'' by analogy with complex geometry,
the surgery is really a symplectic resolution -- but with exceptional
set a Lagrangian three-sphere.  A better name might be a ``Lagrangian
blow-up''. 

\begin{Theorem} \label{smooth}
Every symplectic conifold $(X,\omega_X)$ admits a symplectic smoothing
which contains an embedded Lagrangian $n$-sphere for each node.  In
particular,  
the reverse conifold transition $\tilde{X}$ of any small resolution of
a six-dimensional conifold carries a distinguished symplectic
structure, well-defined up to symplectomorphism.
\end{Theorem}

\begin{Proof}
By Proposition \ref{nodedata} a punctured neighbourhood of each node is
isomorphic to $(\{\sum_iz_i^2=0\}\backslash\{0\},\omega_{\C^{n+1}})$,
which by (\ref{localmodel}) is isomorphic to $T^*S^n\backslash S^n$
with its canonical symplectic structure.
Hence we can replace the node $\{z_i=0\ \,\forall i\}$
by the Lagrangian $n$-sphere and smoothly extend the
form keeping it globally symplectic.

For the uniqueness statement when $n=3$, recall that the set of choices in
performing the $C^\infty$ surgery is
connected. So the resulting smooth manifold is unique,
with two different sets of choices giving two different symplectic
forms on it. These forms are connected by a
family of symplectic forms, coming from connectedness of the space of
gluings, and the cohomology class along the family is constant 
(since the neighbourhood $T^*S^3$ has trivial $H^2$). Moser's theorem
\cite{Mos} then gives the required symplectomorphism.
\end{Proof}

\begin{Remark} \label{limitofforms}
In dimension $n=3$, one could alternatively start with a smooth
six-manifold $X$ with a closed
two-form $\eta$ which is non-degenerate except along a two-sphere $C$,
where it coincides with the appropriate local model (pull-back to
a small resolution of the standard form on the node). This situation
arises naturally as the limit of a path of symplectic forms on $X$
as in the Remarks after Theorem \ref{resol}.
Note that finding
symplectic two-spheres in a simply-connected symplectic
six-manifold $X$ is straightforward; 
they are governed by an $h$-principle \cite{Gro}.
However, finding families of symplectic
forms which degenerate only along such spheres, yielding conifolds,
is more subtle.
\end{Remark}

Secondly, and more deeply, we want to prove the mirror of the results of
\cite{Fr}, \cite{Ti}, which we describe now.

Fix a complex 3-fold $X$ with trivial canonical bundle,
nodal singularities only, and such that small resolutions satisfy the
$\partial \overline{\partial}$-Lemma (these are called 
``cohomologically K{\"a}hler'' by Lu-Tian). This final condition,
\emph{which won't concern us in the mirror situation}, is needed
to be able to use Hodge theory to relate deformations of complex
structure $H^1(TY)$ to 3-cycles $H_3(Y)\cong H^3(Y)$ on a smooth 3-fold
$Y$ with trivial canonical bundle. (Since $K_Y$ is trivial there is an
isomorphism 
$TY\cong\Lambda^{2,0}T^*Y$ and so $H^1(TY)\cong H^{2,1}(Y)\subset H^3(Y)$.)
Following ideas of Clemens \cite{Cl}, Friedman showed
that a necessary condition for the existence of a complex smoothing
of $X$ is that there is a relation in homology between the exceptional
curves $C_i\cong\Pee$ in a small resolution $Y$, of the form
\begin{equation} \label{reln1}
\sum_i\lambda_i[C_i]=0\ \in\ H_2(Y;\Z)
\quad \mathit{with}\ \ \lambda_i\ne0\ \ \mathit{for\ all\ }i.
\end{equation}
(This condition is independent of the choice of small resolution,
as flopping a curve $C_j$ simply reverses the sign of $\lambda_j$
in (\ref{reln1}).) He also showed that (\ref{reln1}) is sufficient
for a first order infinitesimal smoothing of $X$; Tian showed that
this deformation is always unobstructed, i.e. can be extended to
give a genuine smoothing.

Such a \emph{``good''} relation is given by a 3-chain bounding
the $\lambda_iC_i$
in $Y$, which becomes a 3-cycle in $X$ passing through all
of its nodes (and a ``Poincar{\'e} dual'' 3-cycle
in the $C^\infty$ smoothing of $X$, which we may think of as a
vanishing cycle for the nodes).
Intuitively, via the correspondence between 3-cycles
and deformations of complex structure on Calabi-Yau manifolds,
this gives a global deformation of complex structure that restricts
at each node to ($\lambda_i$ times) the unique standard local
smoothing of that node. So for the $\lambda_i$s non-zero the result
is a smooth 3-fold.

The mirror situation for resolutions is perfectly analogous.

\begin{Theorem} \label{resol}
Fix a symplectic 6-manifold $X$ with a collection of $n$ disjoint
embedded Lagrangian 3-spheres $L_i\cong S^3$. 
There is a \emph{``good''} relation (cf. (\ref{reln1}))
\begin{equation} \label{reln2}
\sum_i\lambda_i[L_i]=0\ \in\ H_3(X;\Z) \quad \mathit{with}\ \ \lambda_i\ne0
\ \ \mathit{for\ all\ }i
\end{equation}
iff there is a symplectic structure on one of the $2^n$ choices of
conifold transitions of $X$ in the Lagrangians $L_i$, such that the
resulting $\Pee$s $C_i$ are symplectic.
\end{Theorem}

\begin{Proof}
Via (\ref{localmodel}) we can replace each Lagrangian sphere by a node
and then replace the node by a two-sphere via a small resolution.
This gives a manifold $(Y, \omega)$ where
$\omega$ is globally closed, and degenerate along a collection of
embedded two-spheres $C_i \subset Y$ (i.e. $\omega$ is the form pulled
back from the symplectic conifold). We show that the 4-chain
giving the homology
relation gives rise to a four-cycle $\sigma$ on a resolution (its
$S^3$ boundaries have been collapsed to $S^2$s) with $\sigma\,.\,C_i=
|\lambda_i|>0$. Firstly we give the local model.

In our $(u,v)$-coordinates (\ref{localmodel}) on $T^*S^3$ a
collar neighbourhood of the $S^3$ zero-section $\{\sum u_i^2=1,\ v=0 \}$
is given by the equations defining half of a real line-bundle over
$S^3$ (such a line bundle being necessarily trivial):
$$
\Delta = \{ (u,v) \ | \ v_1 = - \lambda u_2, \ v_2 = \lambda u_1, \ v_3 =
-\lambda u_4, \ v_4 = \lambda u_3; \ \lambda \geq 0 \}.
$$
Using quaternionic multiplication (cf. the Appendix) we can write this
as  $\{ (u,v) \, | \, v=\lambda Iu, \ \lambda \geq 0\}$.  One can check that 
under the diffeomorphism
(\ref{localmodel}) $\Delta$ is exactly the image of the complex surface
$S=\{ z_1 = i z_2, z_3 = iz_4 \}$ inside $\C^4$ which lies inside the
quadric and contains the node; $\lambda$ appears as $\sum\Ree(z_i)^2$.  The
other ``half'' of the line bundle, taking $\lambda \leq 0$ in the
defining equations above, arises from the second complex surface $\{
z_1 = -iz_2, z_3=-iz_4\}$.  The surface $S$
is smooth, as is
its proper transform in either small resolution. Depending on
the resolution chosen (see for instance \cite{LuT}),
this proper transform is either isomorphic
to $S$, intersecting the exceptional $\Pee$ transversally in $+1$, or
is the blow up $\widehat S$ of $S$ at the origin with its
exceptional $\Pee$ coinciding with the exceptional set of the small
resolution; in this case it is easy to see that $\widehat S\,.\,\Pee
=-1$.

By definition our global 4-chain is an element of $H_4(X,\cup_iL_i)$
which maps to $\oplus_i\lambda_i[L_i]\in\oplus_iH_3(L_i,\Z)$ in the
homology exact sequence of the pair. By excision we can replace
$\lambda_i[L_i]\in H_3(L_i,\Z)$ by $\lambda_i\partial(\Delta_i\cap U_i)$
in $H_3(\Delta_i\cap U_i,\Z)\cong H_3(L_i,\Z)$, where $U_i$ is a small 
tubular neighbourhood of $L_i$ in $X$, and $\Delta_i$ is the collar
of $L_i$ defined above in local coordinates. Also by excision, our
global relation gives a four-chain
in $H_4(X,\cup_i(\Delta_i\cap U_i))\cong H_4(X,\cup_iL_i)$ with boundary
$\lambda_i\partial(\Delta_i\cap U_i)$. Adding this to the collars
$(\Delta_i\cap U_i)$ give a 4-chain (homologous to the original
chain) which in local coordinates is exactly $\lambda_i$ times
our collar model around $L_i$. Thus, choosing the right small
resolution (with local intersection of the complex surface with
the exceptional $\Pee$ given by sign$(\lambda_i)$) gives a
4-cycle $\sigma$ with intersection $|\lambda_i|$ with $C_i$.  Write
$\tilde{\sigma}$ for a two-form in the class Poincar{\'e} dual to $\sigma$.

Fix some tubular neighbourhood $U_i$ of each curve $C_i$.
Since $H^2(U_i;\R)\cong\R$ we know that  $\tilde{\sigma}\res{U_i}$ is
cohomologous to $\lambda_i\omega_i$, 
where $\omega_i$ is the standard K{\"a}hler form on the total space
of $\OO(-1)\oplus\OO(-1)$ over $C_i\cong\Pee$ for which $C_i$ has area $1$.

Write $\tilde{\sigma}\res{U_i}=\lambda_i\omega_i+d\phi_i$ on $U_i$, and (via
cut-off functions) pick a
one-form $\phi$ on $Y$ such that $\phi\res{U_i}=\phi_i\ \ \forall i$.
Then replacing $\tilde{\sigma}$ by $\tilde{\sigma}+d\phi$ we may
assume that $\tilde{\sigma}$ 
restricts to $\lambda_i\omega_i$ in a neighbourhood of each $C_i$.

By the compactness of $X \backslash \bigcup U_i$, and the openness
of the non-degeneracy condition, we may choose $N$
sufficiently large that $\Omega=N\omega+\tilde{\sigma}$ is a symplectic form
on $Y\backslash\bigcup_iU_i$. We claim that $N\omega+\tilde{\sigma}$
is in fact a global
symplectic form. As $\omega\res{C_i}=0$ and $\tilde{\sigma}\res{U_i}=
\lambda_i\omega_i$, $\Omega$ is non-degenerate in some smaller neighbourhood
$V_i\subset U_i$ of the $C_i$. The remaining place to check is in
$U_i\backslash V_i$.   Now in general, convex combinations of
symplectic forms are not symplectic, and certainly the forms we have
here are not directly proportional (for instance, $\omega$ is induced
from a form on $T^* S^3$ and is exact, whereas each $\omega_i$ is
non-trivial in cohomology).

However, as Gromov first pointed out, if two symplectic forms $\omega,
\omega_i$ both tame some fixed almost complex structure, then convex
combinations are necessarily symplectic. We are in just such a situation.
Since non-degeneracy is local, using Lemma \ref{nodedata} we may as
well work on a neighbourhood of 
the origin in the standard node $W$ with its standard symplectic
structure. This gives us a complex
structure $J$ -- the standard one on $\C^4$ -- tamed by $\omega$ and
a \emph{holomorphic} resolution $\OO_{\Pee}(-1)^{\oplus2}$ with $\omega_i$
a standard K{\"a}hler form on it. Thus on $U_i\backslash V_i$
both forms tame the same complex structure and so we can
take convex combinations of them.

This completes the ``if'' part of the Proof.  For ``only if'',
certainly if the $C_i$ are symplectic then there is a two-form which
is non-trivial on each.   Any such non-degenerate
form $\omega$ gives a 2-form on $X\backslash\cup_iL_i$,
via the isomorphism (\ref{localmodel}). This fits into the exact
sequence of the pair $(X,X\backslash\cup_iL_i)$ (using the Thom
isomorphism $H^3(X,X\backslash\cup_iL_i)\cong H^0(\cup_iL_i)$)
as follows:
\begin{eqnarray*}
H^2(X\backslash\cup_iL_i)\to \!\!&\!\! \bigoplus_iH^0(L_i)\cong\bigoplus_i
\mathbb R \!\!&\!\! \to H^3(X) \\
\omega \ \qquad \mapsto  \!\!&\!\! \bigoplus_i(\int_{C_i}\omega)
\end{eqnarray*}
Since the third map is (Poincar{\'e} dual to) the inclusion of the
fundamental classes of the $L_i$ into $X$, this gives the required
good relation (\ref{reln2}) $\sum(\int_{C_i}\omega)[L_i]=0$.
\end{Proof}

\begin{Remarks}
As we tend $N \rightarrow
\infty$ above, $\omega + \tilde{\sigma}/N$ is symplectic, and in the limit
degenerates along each of the $C_i$
(cf. the Remark after Theorem
\ref{smooth}), giving us a two-form locally isomorphic to a
pull-back from the node in $\C^4$ to a small resolution.

The proof shows that even if the conditions of the Theorem are not
satisfied, we can always induce a distinguished homotopy class of
almost complex structures on the surgered manifold.  To do this, we
choose some non-degenerate
two-form $\omega_i$ near each $C_i \subset Y$ which extends the
degenerate form $\omega$, and then -- via cut-off functions and the
same tameness and convexity argument -- extend to a non-closed
non-degenerate global two-form. Such forms are in one-to-one
correspondence with homotopy classes of almost complex structures.
\end{Remarks}

Clearly the surgery we have
described does not change the fundamental group of the manifold.
Moreover, if we surger $n$ Lagrangian spheres $L_i$, then we increase
the Euler characteristic by $2n$.  It is then easy to deduce the
following (which is well-known -- see for instance \cite{Cl};
we learnt it from \cite{Gr}):

\begin{Theorem} \label{hom}
If the $n\ L_i$ of Theorem \ref{resol} span an $r$ dimensional
subset of $H_3(X)$, then $b_3(Y)=b_3(X)-2r$, and $b_2(Y)=b_2(X)
+(n-r)$. \hfill$\square$
\end{Theorem}

Intuitively, for every 3-cycle we lose by degenerating the $L_i$, we
lose another by Poincar{\'e} duality ($H_3$ is always even dimensional
and has a symplectic basis). This ``dual'' 3-cycle $L_i'$ intersects
$L_i$
and on the resolution becomes a 3-chain bounding the curve $C_i$; hence
we lose $2r$ lots of $H_3$ and $r$ lots of our new $n$ 2-cycles $C_i$.
(Dually the $n-r$ relations amongst the $[L_i]$ are given by 4-chains
on $X$ that become 4-cycles on $Y$ as their boundaries have been
collapsed; thus $h_4$ also increases by $(n-r)$.)

Since our surgery is an almost complex operation, we can also ask
how the Chern classes of the almost complex structure are affected.
$c_1$ is represented by the zero set
of a transverse section of the canonical bundle $\Lambda^3_\C T^*X$.
We can choose a standard nowhere-vanishing holomorphic section on
$T^*S^3\backslash\{\mathrm{Zero\ section}\}$ in its holomorphic
coordinates (\ref{localmodel}). 
This corresponds to a section on the resolution which extends across
the $\Pee$ by Hartog's theorem; the extension is still
non-vanishing since it is non-zero outside a codimension-two
subvariety $\Pee$.  Thus we can refine $c_1$ of both the smoothing $X$
and resolution $Y$ to lie in
$$
H^2(X,\cup_iL_i)\cong H^2(Y,\cup_iC_i),
$$
mapping to $H^2(X)$ and $H^2(Y)$ respectively. In this sense $c_1$
is preserved by the transition. Thus, in particular,

\begin{Lemma}
(Reverse) conifold transitions preserve the condition $c_1=0$.
\hfill$\square$
\end{Lemma}

In six dimensions, a homotopy class of almost complex
structures is completely determined by the first Chern class.  (The
second Chern class is then determined by the identity $c_2 = (c_1^2 -
p_1)/2$ and the third Chern class is just the Euler characteristic.)
All smooth six-manifolds with $\pi_1 = 0$ and $\mathrm{Tor}(H^*) = 0$
are almost complex, with the almost complex structures indexed by the
integral lift $c_1$ of the Stiefel-Whitney class $w_2$~. The triple of integers
$(c_1^3, c_1c_2, c_3)$ is necessarily of the form $(2\alpha, 24\beta,
2\gamma)$ for suitable integers $\alpha, \beta, \gamma$ and all such
triples are in fact realised by simply connected symplectic manifolds
\cite{Halic}. However, very little is known about the existence
of symplectic manifolds with given Chern classes, in particular with $c_1=0$.

In the next section we will give examples of the symplectic surgeries provided
by Theorems \ref{smooth} and \ref{resol}. It has proved remarkably difficult,
however, to prove definitely that the surgery does not preserve the subclass of K{\"a}hler manifolds. In
the Calabi-Yau context, which, because of mirror symmetry, we would like
to work, there is an obvious obstruction:

\begin{Lemma} \label{notkahler}
Let $X$ be a simply connected symplectic six-manifold with $c_1(X)=0$ and
$b_3(X) = 0$. Then $X$ is not homotopy equivalent to any
K{\"a}hler manifold. \hfill$\square$
\end{Lemma}

The proof is that any simply connected K{\"a}hler manifold with
$c_1=0$ has \emph{holomorphically trivial}
canonical bundle (since by Hodge theory, $H^{0,1}=0$), hence has a
nowhere zero holomorphic three-form $\Omega$. This is
automatically closed and non-zero in $H^3$, since
$\int\Omega\wedge\overline{\Omega}>0$. 
In order to obtain a manifold with $b_3 =0$ from a conifold
transition, one needs a collection of disjoint Lagrangians -- satisfying
a good relation -- spanning $b_3/2$ dimensions in $H_3$ (this is the
maximum possible, by Poincar{\'e} duality).  Smoothly there is no
obstruction to finding such spheres, but the situation in symplectic
geometry is not clear, and in many examples (cf. the next section) the
numerology and geometry seem to conspire precisely to make this
impossible.  Indeed, so many examples ``just fail'' that it is natural
to wonder if there is some obstruction to finding such a collection of
disjoint Lagrangian spheres; it is even natural to wonder if all
symplectic Calabi-Yaus have $b_3 \geq 2$ just like K{\"a}hler Calabi-Yaus.

In this regard, it is worth pointing out that
there are other obstructions to being K{\"a}hler which 
can \emph{never} be violated by conifold transitions.
For instance, Chern-Weil theory implies that for a Calabi-Yau $n$-fold the
$L^2$-norm of the curvature tensor for the Ricci-flat metric is given
by $c_2 \cdot \omega ^{n-2}$, which must therefore be non-negative.

\begin{Prop}
Let $X$ be a symplectic manifold obtained by conifold transitions on a
K{\"a}hler Calabi-Yau 3-fold.  Then $c_2(X) \cdot [\omega_X] > 0.$
\end{Prop}

\noindent\emph{Proof.}
Via the surgery, the second Chern class changes by addition of the $\Pee$s we
introduce in the resolution.  (This follows from symmetry
together with the local computation
of Tian and Yau \cite{TiY} who compute the effect on
$c_2$ of a flop.) By construction our final symplectic
form evaluates positively on these, so the surgery can
only increase the value of $c_2 \cdot [\omega]$. \hfill$\square$ \\

It is possible that every symplectic six-manifold with
$c_1 = 0$ has $c_2 \cdot [\omega] \geq0$ (even if they aren't all
K{\"a}hler).  For Calabi-Yaus with
``large complex structure limit points'', mirror
symmetry gives a topological interpretation to this positivity: the
Calabi-Yau should admit a fibration by Lagrangian tori, and the limiting
locus of critical points of the fibration (as we tend towards the
large complex structure limit point)
should be codimension four and give a distinguished symplectic
cycle representing $c_2$, cf. \cite{SYZ}, \cite{Gr2}. \medskip

\begin{Remarks}
Li and Ruan \cite{LiR} have studied the effect of reverse conifold
transitions on quantum cohomology. The effect of the conifold
transition (\ref{resol}) on $QH^*$ is an interesting open question.
Whereas the reverse conifold transition removes $\Pee$s and their
Gromov-Witten contributions, so the conifold transition removes
Lagrangian $S^3$s and should have a mirror effect on Joyce's
invariant \cite{J}.

Salur \cite{Sa} has considered deformations of (an appropriate
modification of) special
Lagrangian submanifolds in symplectic Calabi-Yaus.  Her results
accordingly apply to manifolds constructed from conifold transitions.
\end{Remarks}

%%%%%%%%%%%%%%%%%%%%%%%%%%%%%%%%%%%%%%%%%%%%%%%%%%%%%%%%%%%%%%%%%%%%%%%%%%

\section{Assorted examples}

In this section, we present various examples of the surgeries. To warm up we shall consider the situation in
two complex dimensions, where the symplectic geometry of ordinary
double points is already interesting. \\

\noindent\textbf{(a).} The
question that shall motivate us is the following:
in the minimal model programme, the (lack of) ampleness 
of the canonical class of a variety plays a fundamental role.
Let $(X,\omega)$ be a symplectic manifold.  When does the canonical
class $K_X \in H^2(X;\Z)$ itself contain symplectic forms?

For four-manifolds the
canonical class is particularly decisive for the global geometry.  A
theorem of Liu \cite{Liu}
asserts that if $X^4$ contains a symplectic surface $C$ with $K_X
\cdot C < 0$ then in fact $X$ is diffeomorphic to a del Pezzo surface,
and indeed following work of Lalonde and McDuff \cite{LM}, the
symplectic form is isotopic to a standard K{\"a}hler form.  In particular,
$X$ is Fano and $-K_X$ (hence also $K_X$) contains K{\"a}hler forms.  
In general, there are some obvious necessary conditions for $K_X$ to
contain symplectic forms: we must certainly have
$(K_X)^2 > 0$. An observation going back to McDuff \cite{McD} 
shows that if $b_+ (X) = 1$ then this necessary condition is in
fact sufficient.
The proof, however, suggests no adaptation
to the general case, and indeed there seems to be nothing special here
about the canonical class: for minimal 4-manifolds with $b_+ = 1$,
every class in $H^2$ of positive square contains symplectic
forms \cite{LL}. This  certainly fails in general (for intriguing examples see 
\cite{Vi}); moreover, Taubes' work in
Seiberg-Witten theory shows that if $b_+ (X) > 1$ then for $K_X$ to
contain symplectic forms, $X$ must contain no smooth $-1$-spheres.

Perhaps surprisingly, the constraints of positive square and
minimality \emph{are} sufficient in the integrable case.  The following was
observed independently by Catanese \cite{Ca} (with a quite different proof).
 
\begin{Prop} \label{symplectic}
If $X$ is a minimal complex surface with $(K_X)^2 > 0$, then the
canonical class $K_X$ contains symplectic forms.  
\end{Prop}

\begin{Proof} 
The result is obvious
for rational surfaces (where, however, the form may not be deformation
equivalent to the usual K{\"a}hler form -- think of $\mathbb{P}^2$).  Using
the classification of surfaces we can therefore assume $X$ is a minimal
complex surface  of general type.  A classical theorem 
\cite{BPV} asserts that for $r \geq 5$ the morphism $|rK_X|: X \rightarrow
\mathbb{P}^N$ defined by a 
multicanonical linear system is an embedding away from
the union of 
all holomorphic $-2$-spheres, which are contracted to
\emph{rational} double point (or A-D-E) singularities.  Suppose first all the
$-2$-spheres are isolated, hence the singularities are nodes.  Then
there are local analytic co-ordinates $z_i$ on projective space near
any node such 
that the image of $X$ is defined by $\{\sum_{i=1}^3 z_i ^2 = 0, z_j =0
\ \forall j>3 \}$.  By an argument of Seidel (\cite{Se}, Lemma 1.7) there is
an isotopy 
of K{\"a}hler forms on $\mathbb{P}^N$, starting with the Fubini-Study
form and compactly supported near each node, which yields a form which
in our analytic co-ordinates is exactly the standard form
$\omega_{\C^N}$ in a small neighbourhood
of the nodes.  (Seidel's argument is local, so we adjust the form on a
ball within the domain of definition of our co-ordinates; the ambient
linear structure of $\mathbb{P}^N$ plays no role.)  Since the isotopy
is through K{\"a}hler forms they are
all non-degenerate on the image of $X$, which is complex.
Pulling back this form 
from projective space to $X$, we get a closed two-form $\eta$ on
$X$ in the cohomology class $K_X$ which is non-degenerate away from
the isolated $-2$-spheres and conforms to the standard model near each
one.  That is,
there is a neighbourhood $U$ of each contracted sphere $C$ such that 
$\eta|_{U \backslash C} \ \cong \ 
\omega_{\C^3}|_{\{\sum z_i ^2 = 0 \} \backslash \{0\} }.$  Using the
fact that

\begin{equation} \label{local}
\left( \big \{ \sum_{j=1}^{3} z_j^2 = 0 \big \} \backslash \{ 0 \}, 
\omega_{\C^3} \right) \ \cong \ \big( T^* S^2 \backslash
S^2, \omega_{can} \big). 
\end{equation}

\noindent we can extend the closed 
form $\eta$ over $U$ as a global symplectic form, making the
two-sphere Lagrangian. 
Since the extended  form vanishes on the two dimensional
homology of $T^* S^2$ -- indeed, the form is exact over $U$ --  it
represents the same cohomology class as $\eta$, as required.

For the general case, consider the total space of a smoothing $\{
f(z_1, z_2, z_3) = 
t \}_{t \in \C}$ of an ADE singularity $\{ f(z_1, z_2, z_3) = 0 \}$,
equipped with the restriction of the standard K{\"a}hler form.  By Seidel's
result on isotopies of K{\"a}hler forms as above, the
symplectic form on $X$ induced from the Fubini-Study form is smoothly
isomorphic to the
standard form in an open neighbourhood of the contracted spheres,
minus the spheres themselves.  
As in \cite{Se}, there is a symplectic parallel transport on the
fibres of the smoothing, which -- restricting to a ray in the base
starting at zero -- shows that the complement of the ADE singularity
in the zero-fibre is symplectomorphic to the complement of the
vanishing cycles in a nearby fibre.  This
general fibre contains a tree of Lagrangian spheres; the contact boundary of
such a domain is always $\omega$-convex \cite{Et}.  Gray's stability
theorem for contact structures \cite{Gra} shows that the boundary of
the neighbourhood of the ADE chain in the central fibre is
contactomorphic to this; the isomorphism of symplectic neighbourhoods
coming from the parallel transport shows it is also $\omega$-convex.
A theorem of
Etnyre \cite{Et} allows one to glue symplectic domains in this
setting, so we can replace a neighbourhood of the
ADE chain by a neighbourhood of the vanishing cycles in a nearby
fibre. (Differentiably this is a trivial surgery - we remove a plumbed
neighbourhood of a tree of spheres and then replace it by a
diffeomorphism of the boundary which is isotopic to the identity.)  We
can therefore 
symplectically extend the two-form $\eta$ pulled back from projective
space with a global symplectic form $\eta'$ for which all the spheres are
Lagrangian, hence has unchanged cohomology class.
\end{Proof}

The canonical class contains K{\"a}hler forms only if
the surface has no holomorphic $-2$-spheres, so this result is a purely
symplectic phenomenon. It would be interesting to use Proposition
\ref{symplectic} as an obstruction to
integrability, for instance on a homotopy K{\"a}hler manifold with 
$\pm K_X$ the only Seiberg-Witten basic classes. \\

\noindent\textbf{(b).} 
We now return to the 3-fold case and give examples of Lagrangian $S^3$s
with good relations (\ref{resol}), therefore admitting symplectic conifold
transitions. As a simple local case, consider the Lagrangian $S^3\subset\Pee\times\C^2$
given by the product of the Hopf map and the unit norm inclusion. Here we have to take the symplectic structure on $\Pee$ given by \emph{minus} that
coming from symplectic reduction of the $S^3\subset\C^2$. Alternatively,
we can compose the Hopf map with the antipodal map, and use the symplectic
structure compatible with the usual complex structure. Taking smaller $S^3$s
in $\C^2$ (and so smaller symplectic forms on $\Pee$), and using the Darboux
theorem to make any K\"ahler surface look locally symplectically like $\C^2$, we get

\begin{Lemma}
Fix a K\"ahler surface $S$, and denote by  $\omega_S$ the pullback of its
K\"ahler form to $\Pee\times S$. By taking $\varepsilon$ sufficiently small,
we can find arbitrarily many disjoint null-homologous Lagrangian three-spheres
in $\Pee\times S$ with symplectic structure $\omega_S^{\,}+\varepsilon
\omega^{\,}_{\Pee}$. \hfill$\square$
\end{Lemma}

(We have only left to show that the spheres are null-homologous: in the original model the sphere bounded the 4-chain $\{\,([v],v)\in\Pee\times\C^2\
:\ 0<|v|<1\,\}$.)

So we can apply Theorem \ref{resol} to these examples to produce, for instance,
three-folds with arbitrarily high $b_2$ and $b_3=0$ which are not obviously
blowups of smooth 3-folds. In the $\Pee\times\mathbb P^2$ case, it should
be possible to show using standard projective 3-fold theory that the degeneration
to a single node and resolution cannot be K\"ahler \cite{Co}, but there are
still some points to check. \\

\noindent\textbf{(c).} 
Though many of the above surgeries are surely 
not realisable within K\"ahler geometry, we have so far been unable to prove it in a particular case. If one could
find a \emph{homotopically} trivial Lagrangian $S^3$ bounding an embedded
$D^4$, then non-K\"ahler manifolds would certainly result.

\begin{Lemma}
The symplectic conifold transition in a Lagrangian $S^3$ bounding an embedded
$D^4$ would violate Hard Lefschetz.
\end{Lemma}

\begin{Proof}
Recall that the Hard Lefschetz theorem for K\"ahler 3-folds \cite{GH} implies in particular that the intersection pairing 
$$
\cap\,[PD(\omega)]: H_4(X) \rightarrow H_2(X)
$$
is an isomorphism.  In
particular, if there is some element $D \in H_4(X)$ for which $D \cap
D' = 0 \in H_2$ for \emph{every} $D' \in H_4$,
then $X$ cannot be K{\"a}hler.  In our case, the class $D$ comes from
the four-ball bounding the Lagrangian sphere.

The single homology relation $[L] = 0$
satisfies the conditions of the surgery Theorem \ref{resol}; following
the proof of that result, the symplectic structure on the resolution $X$
is obtained by deforming with the Poincar{\'e} dual of
a four-cycle which is the lift of the
bounding topological four-ball which is transverse to the $\Pee$ at
one point.  Topologically, this lift is just an embedded four-sphere $D$ inside
the resolution.  For any other four-cycle $D'$, the intersection
product $D \cap D' \in H_2(X)$ is geometrically represented by a
two-cycle lying inside $D$, hence lies in the image of $H_2(D)
\rightarrow H_2(X)$, but this is trivial since $D$ is a sphere.
Hence, the conifold transition of $Z$ along $L$ violates the Hard
Lefschetz theorem and is not homotopy K{\"a}hler.
\end{Proof}

Unfortunately, we do not know of any examples of even homologically
trivial Lagrangian three-spheres in Calabi-Yau 3-folds.  (The known
constructions of such spheres yield submanifolds whose Floer homology
is probably not well-defined, whereas the Floer homology of any pair of
Lagrangian three-spheres in a symplectic Calabi-Yau is well-defined by
work of Fukaya, Oh, Ohta and Ono.) In this direction, however, we give a
new construction of null-homologous Lagrangian 3-spheres in certain fibre products; these are therefore not simple $\Pee\times S$ product examples
as in Example b. \\

\noindent\textbf{(d).} 
The idea of this construction is to represent $S^3$ as a family of tori $T^2$
over an interval, shrinking to circles at the end-points.  Provided the
circles that collapse at the two ends span the first homology of
$T^2$, the closed three-cycle will be topologically a sphere,
represented as a genus one Heegaard splitting.  To obtain Lagrangian
spheres this way is slightly more subtle, but one good case is where
the $T^2$-fibres are themselves Lagrangian submanifolds in the fibres of a four-torus-fibred K\"ahler 3-fold.  Such
three-folds arise naturally from fibre products of elliptic surfaces,
as in the work of Schoen \cite{Sch2}.

Let $\pi_1: S_1 \rightarrow \Pee$ and $\pi_2: S_2 \rightarrow \Pee$ be
elliptic fibrations with smooth generic fibres, only nodal singularities
in fibres and with sections (and hence no multiple fibres).  The fibre
product $S_1 \times_{\Pee} S_2 \rightarrow \Pee$ is smooth except over
points of $\Pee$ which are critical for both $\pi_i$, where at
the points $(\mathrm{Node}, \mathrm{Node})$ there are ordinary double
points, which therefore admit small resolutions.  In particular,
in the generic case, if we start with two rational elliptic surfaces
each of which has 12 singular fibres, and no singular fibres are in
common, the fibre product is a smooth algebraic 3-fold, in fact
Calabi-Yau \cite{Sch2}.  Observe that if the original elliptic
surfaces are equipped with K{\"a}hler forms, and if there are no common
singular fibres, there is a natural K{\"a}hler
form on the fibre product, restricted from $S_1 \times S_2$.

Fix points $a,b \in \Pee$ such that $a \in
\mathrm{Crit}(\pi_1)$ and $b \in \mathrm{Crit}(\pi_2)$, and an arc
$\gamma$ in
$\Pee$ from $a$ to $b$ (disjoint from all other critical values).
Suppose for simplicity there is at most one node in each fibre.
There is a symplectic parallel transport in each $S_i$, which gives
rise to Lagrangian thimbles $\Delta_i$ ($i=1,2$) over the arc $\gamma$
such that  in $S_1$
the thimble $\Delta_1$ has boundary a smooth circle inside $(S_1)_b$ and passes
through the node in the fibre $(S_1)_a$, whilst inside $S_2$ the
thimble $\Delta_2$ contains the node in $(S_2)_b$ and has boundary a circle
inside $(S_2)_a$.  Now consider the fibre-products of the
thimbles (see Figure 1),

$$L \ = \ \Delta_1 \times_\gamma \Delta_2\ =\ (\Delta_1 \times \Delta_2) \cap S_1 \times_{\Pee} S_2 \ \subset
\ S_1 \times S_2.$$

\noindent This is a family of two-tori over $\gamma$, one circle of
which collapses at $a$ and the other at $b$; moreover, the natural
product K{\"a}hler form on $S_1 \times S_2$ restricts trivially to $L$
and hence so does the K{\"a}hler form on the fibre product. A local computation
shows that $L$ is smooth;  this is obvious away from the end-points of
$\gamma$, and over the end-points the result follows from smoothness
of the thimbles at their origins, which is well-known \cite{Se}.  Hence
we have constructed a
Lagrangian three-sphere inside the 3-fold $S_1 \times_{\Pee} S_2$.

\begin{figure}[ht] \vspace{1cm}
\begin{center}
\input{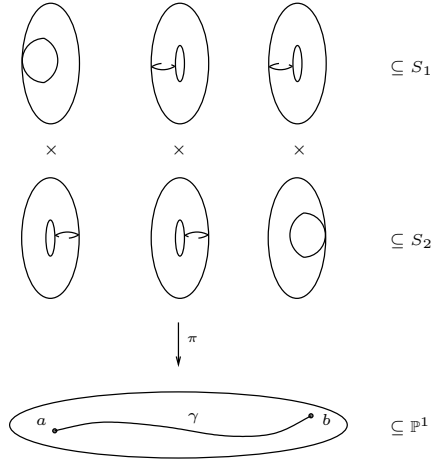}
\end{center}
\caption{Fibred Lagrangian three-spheres} 
\end{figure}

Taking suitable elliptic fibrations, one can again obtain many
examples this way (and often the computation of intersection pairings
can be reduced to counting intersections of arcs inside the base
$\Pee$).  There is an interesting special case, when the vanishing
cycles for the original elliptic surface are homotopically trivial.

\begin{Lemma}
If the elliptic fibrations $S_1$ and $S_2$ have no common
singular fibres and each have a homotopically
trivial vanishing cycle, the (smooth K{\"a}hler) fibre product $Z$ contains
a homologically trivial Lagrangian three-sphere.
\end{Lemma}

\begin{Proof}
In the above construction of the Lagrangian 3-sphere, $\gamma_1$ bounds a
unique disc in each fibre; putting these together gives a 3-chain $D\subset
S_1$, $D^2$-fibred over $\gamma$ except at $a$ where its fibrewise boundary collapses and the $D^2$ becomes an $S^2$ -- the rational component of the singular fibre $(S_1)_a$. Then we have
$$
\partial(D\times_\gamma\Delta_2)\ =\ \Delta_1\times_\gamma\Delta_2\ +\ S^2\times
(\Delta_2)_a\ =\ L\ +\ S^2\times(\Delta_2)_a.
$$
Since the vanishing cycle $(\Delta_2)_a$ was also assumed to be null-homotopic,
its product with the rational component $S^2$ of $(S_1)_a$ is null-homologous;
thus $L$ is also null-homologous.
\end{Proof}

It is straightforward to find suitable elliptic surfaces -- if $\pi:S
\rightarrow \Pee$ is
not relatively minimal (so some fibres contain exceptional $(-1)$-spheres,
i.e. are reducible with one component an $S^2$ meeting only one other
at a single point),
then the vanishing cycles associated to these nodes in the fibres are
homotopically trivial.  Taking a holomorphic automorphism $\phi$ of $\Pee$
that does not fix any of the critical points of $E$, the fibre
product of $\pi$ and $\phi \circ \pi$ contains a Lagrangian sphere as above.
Note that if we blow down the $(-1)$-sphere in $E$, the image of the
thimble is not Lagrangian with respect to any K{\"a}hler form downstairs.

Besides their intrinsic interest, such Lagrangians give interesting
applications of the surgery, since they automatically satisfy a good relation
(\ref{resol}). \\

\noindent\textbf{(e).}
Although we cannot demonstrably
produce non-K{\"a}hler symplectic Calabi-Yaus, we can at least produce
a large collection of symplectic manifolds with $c_1=0$, so many that some
\emph{ought} to be non-K\"ahler; it is frustrating and intriguing that a
year's work has not produced a proof. Here are
some examples that start with the quintic hypersurface in
$\mathbb{P}^5$.  They show that it is possible to apply 
Theorem (\ref{resol}) in concrete cases, and to produce symplectic
manifolds whose Betti numbers are thought not to be
realised by K{\"a}hler Calabi-Yaus. Even if they turned out to be K{\"a}hler,
they would give interesting new examples of rigid (i.e. no complex structure
deformations, equivalent to $h^{2,1}=0$ or $b_3=2$) Calabi-Yaus, many more than are thought to exist.

\begin{Prop}
Symplectic Calabi-Yau manifolds exist with $b_3=2$ and any $2\leq b_2
\leq 25$.
\end{Prop}

\begin{Proof}\,[Sketch]
Consider the hypersurfaces $Q_{\lambda} \subset \mathbb{P}^5$ defined by
\begin{equation} \label{eqn}
\sum_{i=1}^5x_i^5-\lambda\prod_{i=1}^5x_i=0.
\end{equation}
Each has an obvious $(\Z/5)^3$ projective symmetry group $\{(\alpha^{i_1},
\ldots,\alpha^{i_5})\,:\,\sum_ji_j=0$ mod $5\}$.
At $\lambda=5$ this family includes Schoen's quintic \cite{Sch} $Q_5$
with 125 nodes, the $(\Z/5)^3$-orbit of the node at $[1:\ldots:1]$.
For $\lambda\in\R\backslash\{5\},\ Q_\lambda$
is smooth and has a K{\"a}hler structure inherited from $\Pf$. As
symplectic manifolds these $Q_\lambda$ are all isomorphic, and we
call the general such manifold $Q$.  To prove the Proposition, we
exhibited explicit homology relations between subsets of the 125 Lagrangian
vanishing cycles associated to the nodes in $Q_5$.  To find these
homology relations, we begin by constructing enough cycles to span $H_3(Q)$.
These cycles will in fact be unions of real slices: we were heavily inspired by
the calculations  in
Appendix A of the masterpiece \cite{COGP}.

Consider the open chain $\Delta^k$ in $Q_\lambda$ give by taking
$x_1,\ldots,x_4\in(0,\infty)$, and $x_5$ a root (described,
along with $k$, below) of the equation (\ref{eqn}):
\begin{equation} \label{eqn2}
x_5^5=-(x_1^5+\ldots+x_4^5)+\lambda x_1\ldots x_5.
\end{equation}
For $\lambda=0$ we choose the root
$x_5=e^{(2k-1)\pi i/5}(x_1^5+\ldots+x_4^5)^{\frac{1}{5}}$, where the
one-fifth power means the positive real one. For $\lambda\in(0,5)$,
it is shown in Appendix A of \cite{COGP} that
there is a continuous family of choices of $x_5$ compatible with
this one; i.e. the roots of (\ref{eqn2}) do not come together 
as a double root for any given $x_1,\ldots,x_4$ for
$\lambda$ in this range. In fact, differentiating (\ref{eqn2})
shows that the equation for $x_5$ has a double root only when
\begin{equation} \label{branch}
5^5(x_1^5+\cdots+x_4^5)^4=4^4(\lambda)^5(x_1x_2x_3x_4)^5.
\end{equation}
This has no real solutions for $\lambda<5$ (as the geometric mean
of the $x_i^5$ is no larger than the arithmetic mean). For $\lambda>5$
\cite{COGP} show that the $k=0,1$ roots coincide, but we may choose one
(for $x_i$ real) by taking $\lambda\in\R+i\epsilon$ with $\epsilon$ small
and positive (this rules out double roots by (\ref{branch})),
and then tend $\epsilon$ to zero. This defines $\Delta^k$ for
all $\lambda\in[0,\infty)$.

$\Delta^k$ is naturally a 3-simplex with its four boundaries where
one of the $x_i$ ($i\le4$) becomes zero, six edges where two $x_i$
are zero, and four vertices where only one $x_i$ (and $x_5$)
is non-zero. (In fact by scaling the $x_i,\ i=1,\ldots,4$
projectively so that they sum to 1, we get the standard 3-simplex
$\{x_i>0,\ \sum_ix_i=1\}$ in $\R^4$.)

We now describe some similar cells that, along with $\Delta=\Delta^k$,
will make up a simplicial complex in $Q_\lambda\ (\lambda\in(0,5))$,
again for a fixed $k\in\{0,1,2,3,4\}$. Picking $\alpha_i\in
\{0,1,\alpha:=e^{2\pi i/5}\}$, denote by $\Delta_{\alpha_1\alpha_2
\alpha_3\alpha_4}=\Delta_{\alpha_1\alpha_2\alpha_3\alpha_4}^k$ the set
\begin{equation} \label{Delta}
x_i\in\alpha_i(0,\infty),\ i=1,\ldots,4,
\end{equation}
with $x_5$ the $k$th root of (\ref{eqn2}) as above, with argument
tending towards $e^{(2k-1)\pi i/5}$ as $\lambda\to0$.

Thus $\Delta_{1111}$ is our 3-simplex $\Delta=\Delta^k$ above,
and we get $2^4=16$ isomorphic copies by letting the $\alpha_i$
range through $\{1,\alpha\}$. Similarly, allowing one $\alpha_i$
to be zero we get the boundary simplices $x_i=0$ of these
3-simplices; setting two $\alpha_i$s to zero gives the edges,
three the vertices.

We set $C^1$ to be the simplicial complex $\bigcup\Delta_{\alpha_1
\alpha_2\alpha_31}$, and $C^\alpha=\bigcup\Delta_{\alpha_1\alpha_2
\alpha_3\alpha}$; these have common boundary
$C^0=\bigcup\Delta_{\alpha_1\alpha_2\alpha_30}$.
(The four faces of any 3-simplex come from changing an $\alpha_i$ from
an $\alpha$ or a 1 to 0, so 3-simplices meet along a common face if
their three indices differ in only one place.)  Finally we define 
$L^k=C^1\cup_{C^0}C^\alpha$.  We claim this is topologically a three-sphere.
Let $\gamma\subset\C$ be the V-shaped
union of the two halflines $[0,\infty)\cup \alpha[0,\infty)$.
Then $L^k$ is the image of $\gamma_1\times\gamma_2\times\gamma_3
\times\gamma_4\cong\R^4$ in $Q$ under the map $[\gamma_1:\gamma_2:
\gamma_3:\gamma_4:x_5]$ where $x_5$ is the $k$th root of (\ref{eqn2}).
But this map simply divides $\R^4$ by projective rescaling by
$(0,\infty)$ acting radially (and we do not allow the origin
as it doesn't define a projective coordinate), making the image
$\R^4\backslash\{0\}/(0,\infty)=S^3$.

At $\lambda=0$, $\Delta^k$ is in fact \emph{special Lagrangian}
(as the fixed point set of the antiholomorphic involution
$x_i\mapsto\bar{x}_i,\ x_5\mapsto\alpha^{2k-1}\bar{x}_5$ on $Q_0$)
so the $L^k$ are piecewise smooth Lagrangians.  Using the symmetry
group $(\Z/5\Z)^3$, we obtain 625 piecewise smooth 
Lagrangian three-spheres in $Q$.  (Although we do not need the fact,
we remark that any disjoint set can
be \emph{smoothed} as Lagrangians.  To show this either glue in
a near-Lagrangian 
local model along the edges and flow the result to make it Lagrangian
-- the obstruction is just $[\omega]\res L\in H^2(L)$, and so zero --
or tinker explicitly with the defining equations for the $L^k$.)  
We omit the proof of the following computational Lemma:

\begin{Lemma}
The 625 spheres described above span the entire (204 dimensional)
third homology of the quintic. 
\end{Lemma}

This can be proved by explicitly computing the intersection numbers
(indeed the geometric intersections) of
all the cycles with one another: this yields a vast matrix, whose rank
we computed using MAPLE.  Next, one can compute the geometric
intersections between the cycles above and the vanishing cycles for
the nodes in Schoen's quintic $Q_5$.  Given this data, a computer
search also yields 
good relations between disjoint sets of these vanishing cycles.  In
this way, we found disjoint sets of $k$ spheres spanning 101
dimensions in homology, for each $102 \leq k \leq 125$, proving the
Proposition.  (However, computer searches
suggest that no disjoint set of spheres -- drawn from vanishing cycles
or from the piecewise Lagrangian cycles constructed above --  spans a
102-dimensional 
subspace in homology; this suggests we cannot achieve $b_3=0$ this way.)
\end{Proof}

Schoen's nodal quintic $Q_5$ \cite{Sch} -- which has the advantage of being
rigid --  is also the total space of an
abelian surface fibration.  (It is
the relative Jacobian for a certain pencil of genus two curves on
$\Pee \times \Pee$; the base-points of the pencil give a distinguished
model for the compactified Jacobian of the reducible singular fibres.
The vanishing cycle for the node at $[1,1,1,1,1]$ in Schoen's example
is, surprisingly, homologous to the sphere $L^1$ described above, as
can also be checked by computing intersection numbers and using MAPLE.)
However, the singular fibres seem too degenerate to find spheres that lift to
the small resolution in this representation.   Schoen has also
classified which fibre products of rational elliptic surfaces yield
rigid CY 3-folds \cite{Sch2}.  In some
cases, it is again possible to identify the non-zero 3-cycles on the
rigid varieties; in general one can again laboriously check that
none of these rigid fibre products contains a Lagrangian sphere of the
fibred variety 
described above -- the spheres one obtains on the nodal fibre products
always lift to have boundary on the small resolutions. \\

\noindent\textbf{(f).} 
Finally, we mention that one could use Theorem \ref{smooth} to symplectically smooth
a Calabi-Yau (for instance) with node(s) which has no complex smoothing (i.e.
the $\Pee$s in a small resolution satisfy no good relation \cite{Fr}). Again,
we have been unable to find an example where it can be proved that the resulting
manifold has no K\"ahler structure at all. It clearly has no K\"ahler degeneration
back to the nodal variety; perhaps the symplectic degeneration over a disc
can be completed over $\Pee$ to give a route to producing symplectic
non-K\"ahler 8-manifolds.

%%%%%%%%%%%%%%%%%%%%%%%%%%%%%%%%%%%%%%%%%%%%%%%%%%%%%%%%%%%%%%%%%%%%%%%

\section{Appendix: The local model}

In this section we give the local model for smoothings and resolutions of
nodes using homogeneous spaces and the like.
We start with $\R^n$ and work in $\R^n\times(\R^n)^*$ with its natural
symplectic structure. Define
$$
M=\{(a,b)\in\R^n\times(\R^n)^*\,:\,\langle b,a\rangle=0,\ a\ne0\ne b\}/\sim
$$
where $\sim$ is the equivalence relation given by the orbits of the symplectic
$(0,\infty)$-action $(a,b)\mapsto(\lambda a,\lambda^{-1}b)$. Thus $M$ has
a natural symplectic structure, also visible by mapping it as an
(unramified) double cover to the coadjoint orbit
$$
N=\{A\in\mathrm{End}\,\R^n\,:\,\mathrm{tr}\,A=0,\ \mathrm{rank}\,A=1\},
$$
via $(a,b)\mapsto A=b\otimes a$. ($N$ is most obviously an adjoint orbit,
but $\mathfrak{gl}(V)\subset V^*\otimes V$ is self-dual via the trace map.)

To see $M$ as a complex variety we use the flat metric on $\R^n$ and its
dual to find unique representatives in each $(0,\infty)$-orbit with
$|a|^2=|b|^2$, i.e. to write $M$ as
$$
\{(a,b)\in\R^n\times(\R^n)^*\,\backslash\,(0,0)\,:\,\langle b,a\rangle=0,
\ |a|^2-|b|^2=0\}.
$$
If we identify $\R^n\times(\R^n)^*\cong\C^n$ via $z_i=a_i+ib_i$, then
$|a|^2-|b|^2=\sum(a_i^2-b_i^2)=\mathrm{Re}\,\sum z_i^2$ and
$\langle b,a\rangle=\sum a_ib_i=\mathrm{Im}\,\sum z_i^2/2i$, and
we find $M$ is the usual node minus the nodal point:
$$
M=\left\{\sum z_i^2=0\right\}\subset\C^n\backslash\{0\}.
$$

The symplectic and complex structures we have exhibited on $M$
combine to give the K{\"a}hler structure inherited from the above
embedding in $\C^n$; we will see this again via moment maps in the
next section.

Writing $T^*S^{n-1}_a=\{(a,b)\in\R^n\times(\R^n)^*\,:\,
|a|=1,\ \langle b,a\rangle=0\}$, and similarly for
$T^*S^{n-1}_b$ the cotangent bundle of the dual sphere in $(\R^n)^*$,
we have isomorphisms
$$
T^*S^{n-1}_a\backslash S^{n-1}_a\stackrel{\sim}{\longleftarrow}M
\stackrel{\sim}{\longrightarrow}T^*S^{n-1}_b\backslash S^{n-1}_a.
$$
Here the first map is $(a,b)\mapsto(a/|a|\,,\,-|a|b)$ and the second
$(a,b)\mapsto(b/|b|\,,\,-|b|a)$, and the maps are symplectomorphisms
with respect to the canonical symplectic structures
(cf. \ref{localmodel}). Thus we see 
$T^*S^{n-1}$ as a sort of symplectic resolution of the node (rather than
as a smoothing); adding in the (Lagrangian) zero-section $S^{n-1}_a$
or $S^{n-1}_b$ at 
the node resolves it compatibly with the symplectic structure. These
really are the smoothing; in fact in the obvious way they are
$$
T^*S^{n-1}_a\cong\left\{\sum z_i^2=\epsilon\right\} \qquad\mathrm{and}
\qquad T^*S^{n-1}_b\cong\left\{\sum z_i^2=-\epsilon\right\}
$$
where $\epsilon\in(0,\infty)$. Choosing different values of $\epsilon\in
\C^*$ gives symplectic isotopies between these different smoothings,
corresponding to different splittings of $\C^n$ into
$\R^n\oplus(\R^n)^*$, twisting the standard one by $\sqrt\epsilon$.
On looping $\epsilon$ once round $0$ these isotopies have monodromy
the symplectic Dehn twists of \cite{Se}.

Alternatively we can form an \emph{oriented real blow-up}\footnote{Here
``oriented'' means that we divide the normal directions only by
positive real scalars.} $\widehat M$ 
of $M\cup{(0,0)}$ (the double cover of the space of matrices
$\{b\otimes a\,:\,\langle b,a\rangle=0\}$ branched over the zero
matrix), replacing the origin with its link
\begin{equation} \label{symmetricSM}
S(M)=\{(a,b)\,:\,|a|=1=|b|,\ \langle b,a\rangle=0\}.
\end{equation}
(I.e. $\widehat M:=\{(a,b),(x,y)\in M\times S(M)\,:\,|a|x=a,\ |b|y=b\}$.)
But $S(M)$ is the sphere bundle of $T^*S^{n-1}_a$ (and $T^*S^{n-1}_b$)
so that $\widehat M$ is also the oriented blow up of $T^*S^{n-1}$ in its
zero-section $S^{n-1}$, with exceptional set an $S^{n-2}$-bundle
over $S^{n-1}$.

So via the two induced different projections $S(M)\to S^{n-1}_a$ and
$S(M)\to S^{n-1}_b$ we can blow down the blow-up $\widehat M$ via
two different $S^{n-2}$-fibrations to get the two different
symplectic resolutions/smoothings.
\begin{equation} \label{dg1}
\begin{array}{ccccc}
&& \widehat M \\ & \swarrow && \searrow \\
T^*S^{n-1}_a &&&& T^*S^{n-1}_b \\
& \searrow && \swarrow \\ && \mathrm{node}
\end{array}
\end{equation}
(Compare this to the 2 small \emph{complex} resolutions of the 3-fold
node obtained by blowing down the blow-up via two different $\Pee$-fibrations,
but notice in this real case the two resolutions are in fact isomorphic
as there is a whole family of blow downs of $\widehat M$ interpolating
between the above two. Really this situation is more analogous to,
and in fact (the double cover of) a real slice of, the Mukai flop
$T^*\mathbb P^{\,n}\leftrightsquigarrow T^*(\mathbb P^{\,n})^*$ \cite{Le}.)
It is important then that swapping $S^{n-1}_a$ and $S^{n-1}_b$
will \emph{not induce the flop} on the (3-fold) resolution, as they
are isotopic, even though the birational symplectomorphism
$T^*S^3_a\leftrightsquigarrow T^*S^3_b$ of (\ref{dg1}) does
not extend across the zero-sections. We shall see that
the flop actually corresponds to changing orientation on $S^3$.

Alternatively we may blow down our master space $\widehat M$ in a
different way to get the \emph{complex} blow-up of the node.
We now take $(a,b)\in\R^{\,n}\oplus\R^{\,n}$, then via the map
$(a,b)\mapsto a\wedge b\in\Lambda^2\R^{\,n}$ we get a map $S(M)\to
Gr^+$ to the Grassmannian of \emph{oriented} 2-planes in $\R^{\,n}$.
$Gr^+$ is naturally complex (thanks to Simon Donaldson for reminding
us of this) by mapping such an oriented plane to the complex line
$\C.(a+ib)$ in $\C^n$. Extending the metric bilinear
form on $\R^{\,n}$ by linearity
to a quadratic form on $\C^n$ we see as before (from $|a|^2-|b|^2=0=\langle a,b
\rangle$) that the points of $\mathbb P^{n-1}$ (lines in $\C^n$) in the image
of this map are the quadric $\sum z_i^2=0$, and we have exhibited the complex
blow-up of the node -- with exceptional divisor the quadric in
$\mathbb P^{n-1}$ -- 
as a blow down of $\widehat M$, dividing by the complex phase $U(1)$-fibres
rotating $a$ and $b$ around each other in the exceptional locus $S(M)$.

The exceptional loci are all homogeneous spaces, and the maps we have
defined can be described by the following diagram of maps between
their exceptional sets; the varieties themselves are cones or bundles
over these. The left hand side is the symplectic smoothing/resolution
side, birational to the complex resolution right hand side.
\begin{equation} \label{dg2}
\begin{array}{ccccc}
&& \frac{SO(n)}{SO(n-2)}=S(M) \\
& \swarrow && \searrow \\
S^{n-1}=\frac{SO(n)}{SO(n-1)} &&&& \frac{SO(n)}{SO(2)\times SO(n-2)}
=Gr^+ \\
& \searrow && \swarrow \\ && \mathrm{node}
\end{array}
\end{equation}
In the upper half of the diagram, the first arrow is an
$SO(n-1)/SO(n-2)=S^{n-2}$-bundle, with choices 
(e.g. the different smoothings given by $T^*S^{n-1}_a$ and $T^*S^{n-1}_b$)
given by the connected set of choices of embedding $SO(n-2)\into SO(n-1)$.
The second arrow is an $S^1$-fibration, dividing out by complex phase.

\subsection*{The 3-fold case}

In the special case $n=4$ of interest to us, exceptional things happen
on both sides. On the left hand side $S(M) \cong S^3\times S^2$ is a
\emph{trivial} $S^{n-2}$ bundle, but in many different ways. The
easiest way to represent this is via quaternionic geometry;
let $S^2$ be the fixed set of complex structures $J$ on $\R^4$ compatible with
the flat metric and a fixed orientation. Then $S^3\times S^2\cong S(M)$
via $(a,J)\mapsto (a,b=Ja)$. The different splittings of $\C^4$ into
$\R^4\oplus\R^4$ give the different such trivialisations of $S(M)$, and so
isotopies between them, e.g. between
$S^3_a\times S^2$ and $S^3_b\times S^2$. Changing
the orientation, however, changes things more dramatically; first
changing the orientation on $S^3$, and secondly changing $J$. The induced
map on $S^3\times S^2$ is
\begin{equation} \label{themap}
(a,J\,\circ\,)\mapsto (a,\,\circ\ a^{-1}Ja),
\end{equation}
where we are thinking of $a\in\R^4\cong\mathbb H$ as a unit quaternion,
$J$ as a unit \emph{imaginary} quaternion, and we note that post-
(instead of pre-) 
multiplication by unit imaginary quaternions gives the complex structures
of the opposite orientation. Thus the above is the right map, since in both
cases $b$ is $J\circ\,a=a\,\circ\,a^{-1}Ja$.

On the right hand side this splitting of complex structures into two
$S^2$s corresponds to the double cover $SO(4)\to SO(3)
\times SO(3)$, giving two extra projections
$$
Gr^+=\frac{SO(4)}{SO(2)\times SO(2)}\To \frac{SO(3)}{SO(2)}=S^2.
$$
This fits the two \emph{small} complex resolutions of the node
into the diagram (\ref{dg2}).

The easiest description of these maps, and the induced blow down from
$\widehat M$, is to take an oriented plane (or a pair $(a,b)$ of norm one
with $\langle a,b\rangle=0$) to the unique complex structure $J\in S^2$
on $\R^4$ compatible with the metric and orientation, and such that
the plane is preserved by $J$ (equivalently $Ja=b$). Then, as mentioned
above, changing orientation changes $J\circ$ to $\circ\,a^{-1}Ja$
and swaps the $S^2$ factors in the corresponding isomorphism
$Gr^+\cong S^2\times S^2$.   We can also map $Gr^+ \rightarrow
S^2\times S^2$ via taking a plane to the unique $J$ (resp. $J'$) which
preserves the plane and is compatible with the (resp. opposite)
orientation, and thereby see that interchanging $a$ and $b$ interchanges
$J$ and $J'$, since $Ja=b$ whilst $J'b=a$.  (Alternatively we may use
metric and orientation 
to write $\Lambda^2\R^4=\Lambda^+\oplus\Lambda^-$ and identify oriented
planes with sums of unit norm self-dual and anti-self-dual 2-forms
in $S(\Lambda^+)\times S(\Lambda^-)$; changing orientation then swaps
the two factors.)

So we see that \emph{change of orientation on $\R^4$ corresponds to
flopping the small resolution}. Indeed the choice of orientation or
compatible oriented complex structure $J$ on $\R^4$ can be related directly
to the choice of resolution by using $(x_1+iJx_1)$, etc. in the
the choice of factorisation of $x_1^2+x_2^2+x_3^2+x_4^2=0$ that
produces the small resolution as a graph. Here we have chosen a
real slice of $\C^4$ to talk about complex structures $J$ on $\R^4$;
different choices of real structure isotop $J$ but do not change the
$S^2$ it lies in, or the small resolution, just as they isotop the
Lagrangian $S^3$ in the smoothing from, for instance, $S^3_a$ to $S^3_b$.

On the smoothing side this change of orientation
corresponds to changing the orientation of $S^3$, so this
explains again how the signs of the $\lambda_i$ in Theorem \ref{resol}
can be changed (i.e. the orientations on the $L_i$ can be changed)
by flopping.

Finally then, we can describe the surgery that the conifold transition
produces. We glue in an $S^3\times D^3$ to $S^3\times S^2$ to get the
smoothing (and the many isotopic choices of product structure on
$S^3\times S^2$ that we have seen make the various smoothings isotopic),
but glue in a $D^4\times S^2$ to get either resolution. The gluings
for the two resolutions differ via an involution of
$S^3\times S^2$, given by the composition of the map (\ref{themap}),
and change of orientation
on the $S^3$ (unit quaternions) and $S^2$ (unit imaginary quaternions).
Notice (\ref{themap}) is non-trivial even to the $S^2$
as it restricts to the Hopf fibration on $S^3\times\{*\}$; the
diffeomorphism acts by the matrix $\left( \begin{array}{cc} -1 & 0 \\
1 & 1 \end{array} \right)$ on $\pi_3 (S^3 \times S^2) =
\Z(\mathrm{Degree}) \oplus \Z(\mathrm{Hopf})$.  In the symmetric
co-ordinates (\ref{symmetricSM}) on $S(M)$, the map $(a,b) \mapsto
(b,-a)$ corresponds to the monodromy around a half-circle in $\C^*$
(with square inducing the Dehn twist map), whilst the symmetry $(a,b)
\mapsto (b,a)$ is the diffeomorphism coming from the flop.

\subsection*{Moment maps}

Finally we show how our node and its smoothings, and the geometric
structures on them, can be seen in a surprising way via moment maps.
We can form $S^{n-1}$ as the quotient of $\R^{\,n}$ by the dilation
$\R$-action $a\mapsto e^\lambda a$ once we remove the fixed point
$0\in\R^{\,n}$. So we might try to form $T^*S^{n-1}$ as a symplectic
quotient of $T^*\R^{\,n}$ by the induced symplectic $\R$-action
$(a,b)\mapsto(e^\lambda a,e^{-\lambda}b)$ on $\R^{\,n}\times(\R^{\,n})^*$.

The derivative of this action gives a Hamiltonian vector field
$(a,b)\mapsto(a,-b)$ on $\R^{\,n}\times(\R^{\,n})^*$ with Hamiltonian
$h_1=\langle b,a\rangle$. Using the standard complex structure $J$ on
$\R^{\,n}\times(\R^{\,n})^*\cong\C^n$ gives, unusually, another
\emph{Hamiltonian} vector field $(a,b)\mapsto J(a,-b)=(b,a)$ with Hamiltonian
$h_2=\frac12(|b|^2-|a|^2)$.

To form the symplectic quotient we fix a level set $h_1=\epsilon_1$
(which effectively divides by the second $\R$-action induced by $h_2$)
and divide this by the first action of $\R$. Alternatively,
to take this second quotient, we could instead just fix a level set of
$h_2=\epsilon_2$. So, defining the complex moment map $h=2i(h_1+ih_2)
=|a|^2-|b|^2+2i\langle b,a\rangle=\sum(a_j+ib_j)^2=\sum z_j^2$,
and setting it equal to $\epsilon=2i(\epsilon_1+i\epsilon_2)$, we arrive
at the quadric
$$
\sum_j z_j^2=\epsilon, \qquad z_j=a_j+ib_j.
$$
So the node ($\epsilon=0$) and its smoothings all arise in this way,
with their canonical symplectic structures (as quotients) and different
complex structures (depending on the level $\epsilon$ we picked)
both restricted from the ambient K{\"a}hler structure on $\C^4$.

Even more strangely, these two $\R$-actions do not commute, so do not
arise from a holomorphic $\C$-action (which is why there is
\emph{no canonical} complex structure
on the quotient). This is more familiar in the context of
hyperk{\"a}hler quotients, where the complex group action
does not complexify to give a ``quaternionic'' group action, but
moment maps nonetheless exist. In fact our situation is
the real slice of just such a situation \cite{Hi}, \cite{Le},
and this is where our extra Hamiltonian, or moment map, comes from;
the non-commutativity of the two $\R$-actions arises from the
non-commutativity of the quaternions that have featured throughout this
Appendix.

\vskip 1cm

\small \noindent {\tt I.Smith@dpmms.cam.ac.uk} \newline
\noindent Department of Pure Mathematics and Mathematical Statistics, 
University of Cambridge, Wilberforce Road, Cambridge, CB3 0WB. UK. \newline\newline
\noindent {\tt richard.thomas@imperial.ac.uk} \newline
\noindent Department of Mathematics, Imperial College, Huxley
Building, 180 Queen's Gate, London, SW7 2AZ. UK. \newline\newline
\noindent {\tt yau@math.harvard.edu} \newline
\noindent Department of Mathematics, Harvard University, One
Oxford Street, Cambridge MA 02138. USA. \\


\begin{thebibliography}{COGP}

\bibitem[BPV]{BPV}
Barth, W., Peters, C. and Van de Ven, A.
\emph{Compact complex surfaces.} Springer-Verlag, Berlin, 1984.

\bibitem[COGP]{COGP}
Candelas, P., de la Ossa, X., Green, P., and Parkes, L. (1991).
\emph{A pair of Calabi-Yau manifolds as an exactly soluble
superconformal theory.} Nucl. Phys. \textbf{B 359}, 21--74.

\bibitem[Ca]{Ca}
Catanese, F. (2002)
\emph{Symplectic structures of algebraic surfaces and deformation.}
Preprint math.AG/0207254.

\bibitem[Cl]{Cl}
Clemens, H.C. (1983).
\emph{Double solids.}
Advances in Math. \textbf{47}, 107--230. 

\bibitem[Co]{Co}
Corti, A. (2003). Private communication.

\bibitem[Do]{Do}
Donaldson, S.K. (2000).
\emph{Polynomials, Vanishing Cycles and Floer homology.}
In \emph{Mathematics: Frontiers and Perspectives}. AMS.

\bibitem[Et]{Et}
Etnyre, J. (1998) 
\emph{Symplectic convexity in low-dimensional topology.}  Topology
Appl. \textbf{88}, 3--25. 


\bibitem[Fr]{Fr}
Friedman, R. (1986).
\emph{Simultaneous resolution of threefold double points.}
Math. Ann. \textbf{274}, 671--689. 

\bibitem[Gra]{Gra}
Gray, J.W. (1959).
\emph{Some global properties of contact structures.}
Ann. of Math. \textbf{69}, 421-50.

\bibitem[GH]{GH}
Griffiths, P. and Harris, J.
\emph{Principles of algebraic geometry.} Wiley, New York, 1978.

\bibitem[Gr]{Gr}
Gross M. (1995).
\emph{Blabber about black holes.} Unpublished notes. 

\bibitem[Gr2]{Gr2}

Gross, M. (1998)
\emph{Special Lagrangian fibrations I: Topology.}
Integrable systems and algebraic geometry (Kyoto, 1997), 156--93, World
Scientific. 

\bibitem[Gro]{Gro}
Gromov, M. 
\emph{Partial differential relations.}  Springer, Berlin, 1989.

\bibitem[Ha]{Halic}
Halic, M. (1999).
\emph{On the geography of symplectic 6-manifolds.}  Manuscripta
Math. \textbf{99}, 371--81. 

\bibitem[H]{H}
Hatcher, A. E. (1983). 
\emph{A proof of a Smale conjecture, $\mathrm{Diff}(S^3)\simeq\mathrm{O}(4)$.} 
Ann. of Math. \textbf{117}, 553--607.  

\bibitem[Hi]{Hi}
Hitchin, N. J. (1992).
\emph{Hyper-K{\"a}hler manifolds}. S{\'e}minaire Bourbaki, Vol. 1991/92.
Ast{\'e}risque \textbf{206}, 137--166. 

\bibitem[J]{J} Joyce, D. D.
\emph{On counting special Lagrangian homology 3-spheres.}
Preprint hep-th/9907013. 

\bibitem[LM]{LM}
Lalonde, F. and McDuff, D. (1996).
\emph{J-curves and the classification of rational and ruled symplectic
4-manifolds.}  Contact and symplectic geometry, 7-47, CUP, 1996.
 
\bibitem[Le]{Le}
Leung, N.-C. (2001).
\emph{Lagrangian submanifolds in Hyperk{\"a}hler manifolds, Legendre
transformation.} Preprint math.SG/0110330. 

\bibitem[LiR]{LiR}
Li, A.-M. and Ruan, Y. (2001).
\emph{Symplectic surgery and Gromov-Witten invariants of Calabi-Yau
3-folds I.} Preprint math.AG/9803036. 

\bibitem[LL]{LL}
Li, T.-J. and Liu, A.-K. (2001).
\emph{Uniqueness of symplectic canonical class, surface cone and
  symplectic cone of 4-manifolds with $b^+=1$.}
J. Diff. Geom. \textbf{58}, 331-70.


\bibitem[Liu]{Liu}
Liu, A.-K. (1996)
\emph{Some new applications of the general wall-crossing formula.}
Math. Res. Lett. \textbf{3}, 569-85.
  

\bibitem[LuT]{LuT}
Lu, P. and Tian, G. (1996).
\emph{The complex structures on connected sums of $S^3\times S^3$.}
Manifolds and geometry (Pisa, 1993), 284--293, Sympos. Math., XXXVI,
CUP.

\bibitem[McD]{McD}
McDuff, D. (1998).
\emph{From symplectic deformation to isotopy.}
Topics in symplectic $4$-manifolds (Irvine, 1996), International Press.

\bibitem[McD-S]{McD-S}
McDuff, D. and Salamon, D. 
\emph{Introduction to symplectic topology.} 2nd ed, OUP, 1998.


%\bibitem[Mc]{Mc}
%McRae, A. (1994).
%\emph{Darboux theorems for pairs of submanifolds.}
%Ph.D. Thesis, SUNY Stony Brook.

\bibitem[MoSz]{MoSz}
Morgan, J. and Szabo, Z. (1997).
\emph{Homotopy $K3$ surfaces and mod two Seiberg-Witten invariants.}
Math. Res. Lett. \textbf{4}, 17--21. 

\bibitem[Mo]{Mo}
Morrison, D. (1999).
\emph{Through the looking glass.}
Mirror symmetry III, 263--277, AMS/IP Stud. Adv. Math., 10,
and alg-geom/9705028. 

\bibitem[Mos]{Mos}
Moser, J.K. (1965).
\emph{On the volume elements on manifolds.}
Trans. Amer. Math. Soc. \textbf{120}, 280--96. 

\bibitem[Re]{Re}
Reid, M. (1987).
\emph{The moduli space of $3$-folds with $K=0$ may nevertheless be
irreducible.} Math. Ann. \textbf{278}, 329--334. 

\bibitem[Sa]{Sa}
Salur, S. (2000).
\emph{Deformations of Special Lagrangian submanifolds.}
Comm. Contemp. Math. \textbf{2}, 365--72. 

\bibitem[Sch]{Sch}
Schoen, C. (1986).
\emph{On the geometry of a special determinantal hypersurface
associated to the Mumford-Horrocks vector bundle.}
J. Reine Angew. Math. \textbf{364}, 85--111. 

\bibitem[Sch2]{Sch2}
Schoen, C. (1988).
\emph{On fiber products of rational elliptic surfaces with section.}
Math. Z. \textbf{197}, 177--199. 

\bibitem[Se]{Se}
Seidel, P. (2001).
\emph{A long exact sequence for symplectic Floer cohomology.}
Preprint math.SG/0105186. 

\bibitem[SYZ]{SYZ}
Strominger, A., Yau, S.-T. and Zaslow, E. (1996).
\emph{Mirror symmetry is $T$-duality.}
Nuclear Phys. B. \textbf{479}, 243--259. 

\bibitem[Ti]{Ti}
Tian, G. (1992).
\emph{Smoothing $3$-folds with trivial canonical bundle and ordinary
double points.} Essays on mirror manifolds, 458--479, Internat. Press,
Hong Kong. 

\bibitem[TiY]{TiY}
Tian, G. and Yau, S.-T. (1987).
\emph{Three-dimensional algebraic manifolds with $c_1=0$ and
  $\chi=-6$.}
Mathematical Aspects of String Theory, 574-628.  World Scientific.

\bibitem[Vi]{Vi}
Vidussi, S. (2002)
\emph{Norms on the cohomology of a 3-manifold and SW theory.}
Preprint math.GT/0204211.

\bibitem[We]{We}
Weinstein, A. (1971).
\emph{Symplectic manifolds and their Lagrangian submanifolds.}
Advances in Math. \textbf{6}, 329--346.

\end{thebibliography}
\end{document}